\documentclass[11pt,leqno]{article}
\usepackage{ graphicx, amsfonts, amsthm, amsxtra, verbatim, multicol}
\usepackage[mathscr]{euscript}
\begin{document}

\begin{center}
{\LARGE\bf Moduli of polarized Hodge structures
}
\\
\vspace{.25in} {\large {\sc Hossein Movasati}} \\
Instituto de Matem\'atica Pura e Aplicada, IMPA, \\
Estrada Dona Castorina, 110,\\
22460-320, Rio de Janeiro, RJ, Brazil, \\
{\tt hossein@impa.br \\
http://www.impa.br/$\sim$hossein/}
\end{center}

\newtheorem{theo}{Theorem}
\newtheorem{exam}{Example}
\newtheorem{coro}{Corollary}
\newtheorem{defi}{Definition}
\newtheorem{prob}{Problem}
\newtheorem{lemm}{Lemma}
\newtheorem{prop}{Proposition}
\newtheorem{rem}{Remark}
\newtheorem{conj}{Conjecture}
\def\K{{\mathbb K}}                  
\def\lieg{{\mathfrak g}}   
\def\podu{{\sf pd}}   
\def\per{{\sf pm}}      
\def\perr{{\sf q}}        
\def\perdo{{\cal K}}   
\def\sfl{{\mathrm F}} 
\def\sp{{\mathbb S}}  
 
\newcommand\diff[1]{\frac{d #1}{dz}} 
\def\End{{\rm End}}              
\def\hol{{\rm Hol}}
\def\sing{{\rm Sing}}            
\def\spec{{\rm Spec}}            
\def\cha{{\rm char}}             
\def\Gal{{\rm Gal}}              
\def\jacob{{\rm jacob}}          
\def\tjurina{{\rm tjurina}}      
\newcommand\Pn[1]{\mathbb{P}^{#1}}   
\def\Ff{\mathbb{F}}                  
\def\Z{\mathbb{Z}}                   
\def\Gm{\mathbb{G}_m}                 
\def\Q{\mathbb{Q}}                   
\def\C{\mathbb{C}}                   
\def\O{{\cal O}}                     
\def\as{\mathbb{U}}                  
\def\ring{{\mathsf R}}                         
\def\R{\mathbb{R}}                   
\def\N{\mathbb{N}}                   
\def\A{\mathbb{A}}                   
\def\uhp{{\mathbb H}}                
\newcommand\ep[1]{e^{\frac{2\pi i}{#1}}}
\newcommand\HH[2]{H^{#2}(#1)}        
\def\Mat{{\rm Mat}}              
\newcommand{\mat}[4]{
     \begin{pmatrix}
            #1 & #2 \\
            #3 & #4
       \end{pmatrix}
    }                                
\newcommand{\matt}[2]{
     \begin{pmatrix}                 
            #1   \\
            #2
       \end{pmatrix}
    }
\def\ker{{\rm ker}}              
\def\cl{{\rm cl}}                
\def\dR{{\rm dR}}                

\def\hc{{\mathsf H}}                 
\def\Hb{{\cal H}}                    
\def\GL{{\rm GL}}                
\def\pese{{\sf P}}                  
\def\pedo{{\cal  P}}                  
\def\PP{\tilde{\cal P}}              
\def\cm {{\cal C}}                   
\def\K{{\mathbb K}}                  
\def\k{{\mathsf k}}                  
\def\F{{\cal F}}                     
\def\M{{\cal M}}
\def\RR{{\cal R}}
\newcommand\Hi[1]{\mathbb{P}^{#1}_\infty}
\def\pt{\mathbb{C}[t]}               
\def\W{{\cal W}}                     
\def\gr{{\rm Gr}}                
\def\Im{{\rm Im}}                
\def\Re{{\rm Re}}                
\def\depth{{\rm depth}}
\newcommand\SL[2]{{\rm SL}(#1, #2)}    
\newcommand\PSL[2]{{\rm PSL}(#1, #2)}  
\def\Resi{{\rm Resi}}              

\def\L{{\cal L}}                     
\def\Aut{{\rm Aut}}              
\def\any{R}                          
\newcommand\ovl[1]{\overline{#1}}    

\def\T{{\cal T }}                    
\def\tr{{\mathsf t}}                 
\newcommand\mf[2]{{M}^{#1}_{#2}}     
\newcommand\bn[2]{\binom{#1}{#2}}    
\def\ja{{\rm j}}                 
\def\Sc{\mathsf{S}}                  
\newcommand\es[1]{g_{#1}}            
\newcommand\V{{\mathsf V}}           
\newcommand\WW{{\mathsf W}}          
\newcommand\Ss{{\cal O}}             
\def\rank{{\rm rank}}                
\def\Dif{{\cal D}}                   
\def\gcd{{\rm gcd}}                  
\def\zedi{{\rm ZD}}                  
\def\BM{{\mathsf H}}                 
\def\plf{{\sf pl}}                             
\def\sgn{{\rm sgn}}                      
\def\diag{{\rm diag}}                   
\def\hodge{{\rm Hodge}}
\def\HF{{\sf F}}                                
\def\WF{{\sf W}}                               
\def\HV{{\sf HV}}                                
\def\pol{{\rm pole}}                               
\def\bafi{{\sf r}}
\def\codim{{\rm codim}}                               
\def\id{{\rm id}}                               
\def\gms{{\sf M}}                           
\def\Iso{{\rm Iso}}                           
\def\Ra{\mathrm{Ra}}
\def\nf{g_2}                         
\begin{abstract}
Around  1970 Griffiths introduced the moduli of polarized Hodge
structures/the period domain $D$ and described
a dream to enlarge $D$ to a moduli space of degenerating polarized
Hodge structures. Since in general $D$ is not a Hermitian symmetric domain,
he asked for the existence of a certain automorphic cohomology theory for
$D$, generalizing the usual notion of automorphic forms on symmetric
Hermitian domains. Since then there have been many efforts in the first
part of Griffith's dream  but the second part still lives in darkness.
The objective of the present text is two-folded. First, we give an exposition of the subject.
Second, we give another formulation of the Griffiths problem, based on the classical
Weierstrass uniformization theorem. 

\end{abstract}
\section{Introduction}
I got acquainted with the subject of present  paper, when I was looking for works on 
Abelian integrals in algebraic geometry. My initial aim was to collect some information
and then try to apply them in the study of Abelian integrals in 
differential equations/holomorphic foliations (see \cite{mov0, mov}).
In a canonical way, I found many articles of P. A. Griffiths  around 1970 on periods of
projective manifolds and in particular, his
survey article \cite{gr70}  which contains almost all of his ideas on the subject. 
After a long period of investigation, I observed that I had gone far from my initial aim. But it
was worthfull because, in my opinion, this is a very beautiful piece of mathematics which can
attract  many young mathematicians. Later, I got the idea of modular foliations
(see \cite{ho06-1,ho06-3}) based on the generalized Griffiths domain and did not 
touch more the problems posed by Griffiths.

The conference 
"CIMPA-UNESCO-IPM School on Recent Topics in Geometric Analysis, 2006" 
was a nice occasion for me
to start writing what I had collected so far on the subject. The present text is the main body of
my talks at the mentioned conference and it is mainly
expository.  Its objective is to introduce the reader with the literature on the moduli of polarized 
Hodge
structures and compactification problems after Satake, Baily and Borel. I have tried to introduce
the subject from the point of view of the classical Weierstrass uniformization theorem. This
seems to me the proper way of the realization of Griffiths dream
 on automorphic type functions for the moduli of polarized Hodge structures. 
Let us explain the contents of this text.


In \S\ref{sec1} we sketch the objective of this text in the case of elliptic curves.
This yields to the classical theory of Eisenstein series and modular forms.
In \S\ref{sec2} we recall the Hodge structures on the de Rham cohomologies  of projective
manifolds  and the associated polarizations. In \S \ref{sec3} we  
construct the classifying space of polarized Hodge structures $D$,
called the Griffiths domain, and the action of an arithmetic group $\Gamma_\Z$ on $D$ from the left.
In \S \ref{sec4} we recall Ehresmann's fibration theorem and
then the fact that the period maps  form coefficient spaces to $\Gamma_\Z\backslash D$
satisfy the so called Griffiths transversality.
In \S \ref{sec5} we state the Baily-Borel theorem
 on the unique algebraic structure of 
quotients of symmetric Hermitian domains by discrete arithmetic groups.  Since except in few
cases $D$ is not  a Hermitian symmetric domain,  one cannot apply this theorem to $D$. 
We will mention such few cases which give  origin to the notion of 
 Shimura varieties in algebraic
geometry.  
Finally, in \S\ref{sec6} we give a new formulation of the automorphic type functions corresponding
to families of hypersurfaces.  

I am currently working on the text \cite{ho06-1} in which 
an analytic variety $P$ over the Griffiths domain is constructed and  
many notions of  the present text, like the action of an arithmetic group, period
map, Griffiths transversality and so on  are extended to $P$.

 \section{Elliptic curves}
 \label{sec1}
 In this section we sketch the objective of the present text in the case of polarized Hodge structures
 arising from elliptic curves. The reader who is not familiar with Hodge structures 
 is recommended to read first this section and then the following sections.
 \subsection{Elliptic integrals and elliptic curves}
 Elliptic integrals
 $$
 \int_{I}\frac{dx}{\sqrt{4x^3-t_2x-t_3}}, \ t_1,t_2\in\C, \Delta:=t_2^3-27t_3^2\not =0,
 $$
 where $I$ is some path in $\C$ with
 end points in the roots of $4x^3-t_2x-t_3=0$ or infinity, 
 can be written, up to some algebraic constants, as
 $\int_{\delta}\frac{dx}{y}$, where $\delta\in H_1(E_t,\Z)$ and $E_t$ is an elliptic curve
 in the 
 Weierstrass family of elliptic curves
\begin{equation}
\label{family}
E_t: y^2-4x^3+t_2x+t_3=0, \ t=(t_1,t_2)\in \as_0:=\C^2.
\end{equation} 
A parameter $t$ with $\Delta(t)=0$ corresponds to the singular $E_t$. 
 In fact, after adding the point at infinity to $E_t$ it
 turns to be a compact elliptic curve and by $E_t$ we mean the compact one.
  The de-Rham cohomology (with complex valued 
 differential forms) of $E_t$ is a two dimensional $\C$-vector space generated by 
 the classes $\omega$ and $\bar \omega$ of the differential forms 
 $\frac{dx}{y}|_{E_t}$ and $\frac{d\bar x}{\bar y}|_{E_t}$, respectively, in $H_\dR^1(E_t)$.
 

 \subsection{Polarized Hodge  structures}
 We have the Hodge decomposition 
 $$
 H_\dR^1(E_t):=H^{10}\oplus H^{01}, H^{01}=\overline {H^{10}},
 $$ 
 where $H^{10}$ is the one dimensional $\C$-vector space generated by $\omega$. 
 From another side we have 
 $$
 H^1(E_t,\Z)\otimes \C=H^1(E_t,\C)\cong H_\dR^1(E_t),
 $$
 where $H^1(E_t,\Z)$ is the dual of the $\Z$-module $H_1(E_t,\Z)$. In $H^1(E_t,\Z)$ we have
 the non-degenerated bilinear  map 
 $$
 H^1(E_t,\Z)\times H^1(E_t,\Z)\rightarrow \Z, \psi(a,b)=\int_{E_t}a\wedge b
 \ \ \hbox{up to a constant} 
 $$
  which
 is dual (by Poincar\'e duality, see \cite{gri} p. 59 and \S \ref{intersectionform})
 to the intersection mapping in 
 $H_1(E_t,\Z)$. 
 It is also called a polarization. 
 It satisfies $-\sqrt{-1}
 \psi(\omega,\bar\omega)>0$ which is equivalent to:
 \begin{equation}
 \label{osakareject}
 \Im(z)>0,\ z:=\frac{\int_{\delta_1}\omega}{\int_{\delta_2}\omega},
 \end{equation}
where $\delta_i,\ i=1,2$ is a basis of the $\Z$-module $H_1(E_t,\Z)$. Usually
one takes the symplectic basis of $H_1(E_t,\Z)$, i.e. $\langle \delta_1,\delta_2\rangle =1$, 
and so in this basis the intersection matrix in $H_1(E_t,\Z)$ is 
$$
\Psi:=\mat{0}{1}{-1}{0}.
$$ 
Let $\check \delta_i\in H^1(E_t,\Z), \ i=1,2$ be the dual of $\delta_i$. We have
 $$
 \omega=(\int_{\delta_1}\omega)\check\delta_1+(\int_{\delta_2}\omega)\check\delta_2
 $$
 and so elliptic integrals
 are encoded in an abstract structure consisting of:
 A $\Z$-module $V_\Z$ of rank two  and a polarization $\psi_\Z:V_\Z\times V_\Z\rightarrow \Z$,
 which is a non-degenerated anti-symmetric bilinear map,  and a Hodge structure 
 $$
 V_\C:=V_\Z\otimes \C=H^{10}\oplus H^{01}, \ 
 H^{01}=\overline {H^{10}}, \ \dim H^{01}=\dim H^{10}=1,
 $$
 $$ 
 -\sqrt{-1}\psi(v,\bar v)>0,\ \forall v\in H^{10}. 
 $$
 One call this a polarized Hodge structure of type $\Phi:=
 (m, h^{10}, h^{01}, \Psi)$, where $m=1$ and  
 $h^{10}=h^{01}=1$.

 \subsection{Moduli of polarized Hodge structures}
 \label{moduli}
 Next, we want to construct the moduli of polarized Hodge structures of
 type $\Phi$.  
 Fix a $\Z$-module of rank $2$ and a polarization $\psi_\Z$ on $V_\Z$. 
  Let $D$ be the space of all polarized
 Hodge structures of type $\Phi$. It is in fact isomorphic to 
 the Poincar\'e upper half plane
 $$
 \uhp:=\{x+iy\in\C \mid  y>0\}.
 $$
 The isomorphism $\uhp\rightarrow D$ is constructed in the following
 way: To $z\in \uhp$ we associate the Hodge structure in which $H^{10}$
 is generated by $z\check \delta_1+\check \delta_2$, where $\{\check \delta_1,\check \delta_2\}$ 
 is a basis of $V_\Z$ with $\psi_\Z(\check \delta_1,\check \delta_2)=1$. 
 The group 
 $$
 \SL 2\Z:=\{A\in \Mat(2,\Z)\mid \det(A)=1\}
 $$
 acts on $\uhp$ by $\mat{a}{b}{c}{d}z:=\frac{az+b}{cz+d}$ and the Hodge structures associated
 to $z$ and $Az,\ A\in\SL 2\Z$ are isomorphic. It is not difficult to see that $\Gamma_\Z\backslash 
 \uhp$ is the moduli of polarized Hodge structures of type $\Phi$.

 \subsection{Period map}
 \label{periodmap}
 The multiplicative group $\C^*:=\C\backslash \{0\}$ acts in $\C^2$ in the following 
 way:
 $$
 \lambda\cdot (x,y)=(\lambda^2x,\lambda^3y),\ \lambda\in \C^*, (x,y)\in \C^2.
 $$
 This induces an action of $\C^*$ on $\as_0$:
 $$
 \lambda\bullet (t_2,t_3)=(\lambda^4t_2,\lambda^6t_3), \ \lambda\in\C^*,(t_2,t_3)\in\as_0
 $$
 and an isomorphism $E_{\lambda\bullet t}\cong E_t$ of elliptic curves in the family 
 (\ref{family}). Therefore, $E_t,\ t\in \Pn {2,3}:=\C^*\backslash \as_0$ may be non-isomorphic elliptic curves.
 Note that $\Delta=0$ induces a one point set $c$ in $\Pn {2,3}$.
 We have the period map form $\Pn {2,3}\backslash \{c\}$ to the moduli of polarized Hodge structures of type
 $\Phi$. Composing this with the isomorphism obtained in \S \ref{moduli} and extending it
 to $c$ we have:
 \begin{equation}
 \label{12.5.06}
 \per: \Pn {2,3}\rightarrow (\Gamma_\Z\backslash \uhp)\cup \{\infty\}, \ 
 \per(t)=[z], \ z:=\frac{\int_{\delta_1}\omega_1}{\int_{\delta_2}\omega}, \ \per (c)=\infty,
 \end{equation}
 which we call it again the period mapping. The different choice of the cycles $\delta_1,\delta_2$
 with $\langle \delta_1,\delta_2\rangle=1$ will lead to the action of $\SL 2\Z$ on $z$ and so the above map is well-defined.
 The point $\infty$ can be thought of the point obtained by the action of $\SL 2\Z$ on $\Q$.
 It is called a cusp point.  
\subsection{Gauss-Manin connection}
\label{gm}
For practical purposes it is useful to redefine the period
map in the following way: Let  $\L$ be the set of lattices $\Z\omega_1+\Z\omega_2,\ \frac{\omega_1}{\omega_2}\in \uhp$. The new period map is
\begin{equation}
\label{12.5}
\per: \as_0\backslash \{\Delta=0\}
\rightarrow \L,\ \per(t)=\frac{1}{\sqrt{2\pi i}}\int_{H_1(E_t,\Z)}\omega.
\end{equation}
The previous period map is the quotient of the new one. The derivative
of the period map can be calculated using the following argument: 
For a cycle $\delta\in H_1(E_t, \Z)$ define
$$
\eta_1= \int_{\delta} \frac{dx}{y}, \ \eta_2= \int_{\delta} \frac{xdx}{y}
$$
which are multi-valued holomorphic functions in $\as_0\backslash\{\Delta=0\}$.
Then $\eta_1, \eta_2$ satisfy the following
Picard-Fuchs system
\begin{equation}
\left(%
\begin{array}{c}
  d\eta_1 \\
  d\eta_2 \\
\end{array}%
\right) =\frac{1}{\Delta}
\left(%
\begin{array}{cc}
  -\frac{d\Delta}{12} & -\frac{3 \delta}{2} \\
  - \frac{t_2 \delta}{8} & \frac{d\Delta}{12}\\
\end{array}%
\right)
\left(%
\begin{array}{c}
  \eta_1 \\
  \eta_2 \\
\end{array}%
\right), \ \delta=3t_3dt_2-2t_2dt_3.
\end{equation}
(see for instance \cite{sas74}). The $2\times 2$ matrix in the above
equality is called the Gauss-Manin connection matrix of the family (\ref{family})
with respect to the differential forms $\frac{dx}{y},\ \frac{xdx}{y}$.
The algorithms which calculate the Gauss-Manin connection can be implemented
in any software for commutative algebra
(see \cite{hos005, ho06-1}). 
 \subsection{Eisenstein series and full modular forms}
 For a given elliptic curve in the Weierstrass family 
  we constructed the associated Hodge structure. Now,
 is it possible to construct the inverse map? In fact, it turns 
 out that the period map (\ref{12.5.06}) is a biholomorphism 
  whose inverse is
 given by the $\SL 2\Z$ invariant $j$-function:
 $$
j(z):=q^{-1}+744+
196884 q+21493760q^2+\cdots, \ q=e^{2\pi iz}.
$$
Let us sketch the proof. The period map in (\ref{12.5}) is a local biholomorphic map.
This can be check by the formula of its derivative.
It satisfies also
$$
\per(\lambda\bullet t)=\lambda\per(t).
$$
Therefore, the induced map $\Pn {2,3}\backslash\{c\}\rightarrow \SL 2\Z\backslash \uhp$ is  a local
biholomorphism. Since $\Pn {2,3}\backslash \{c\}\cong \C$, it must be a biholomorphism. 

The  inverse of the period map (\ref{12.5}) composed with the canonical map 
$\{(\omega_1,\omega_2)\mid \frac{\omega_1}{\omega_2}\in \uhp\}\rightarrow \L$ 
is given by $(g_4,g_6):=(60E_4, -140E_6)$, 
where $E_i$ is the Eisenstein
series of weight $i$:
$$
E_i(\Lambda):=\sum_{0\not =a\in \Lambda}\frac{1}{a^{i}}=\sum_{0\not =(n,m)\in\Z^2}
\frac{1}{(m\omega_1+n\omega_2)^i}.
$$
This follows from Weierstrass uniformization theorem (see for instance \cite{si94}).
We have
$$
j(z)=\frac{g_4^3(z)}{g_4^3(z)-27g_6^2(z)}.
$$
A holomorphic function $g$ in $\L$ is called a full modular form of weight $i$, $i$ an even positive integer,  if 
$$
f(k\cdot \Lambda)=k^{i}f(\Lambda),\ \lambda\in \C^*,\Lambda\in \L
$$
and the pull-back $\tilde f$ of $f$ by the canonical map $\uhp\rightarrow \L$ has a finite
growth at infinity, i.e. $\lim_{\Im(z)\rightarrow +\infty} \tilde f(z)=a<\infty$.
The function $g_i,\ i=4,6$ is a full modular form of weight $i$. The fact that
the period map (\ref{12.5}) is a  biholomorphism 
implies that 
every full modular form can be written in a unique way  as a polynomial in $g_4$ and $g_6$. 
It is easy to see that there is no full modular form of odd weight.

The classical definition of a full modular form is as follows: A holomorphic 
function $\tilde f$ on $\uhp$ is called a modular form of weight $i$ if it has a finite growth at
infinity and
$$
(cz+d)^{-i}f(Az)=f(z),\ \forall z\in \uhp,\ A\in \SL 2\Z.
$$
The map $f\mapsto \tilde f$ constructed in the previous paragraph is a bijection
between the two notions. 

One can interpret the modular forms of weight $i$ (in the second sense) as the sections 
$K^{\frac{i}{2}}$,
where $K$ is the canonical  bundle of $\SL 2\Z\backslash \uhp$(see \S \ref{autfun}).
The literature of modular forms and their arithmetic properties is huge. For instance,
the Fermat last theorem is proved to be equivalent to Shimura-Taniyama conjecture
which asks for the existence of certain modular forms. The moonshine conjecture
interprets the coefficients of the $j$-function in terms of the representation
of  the Monster group.

 \section{Hodge theory of projective manifolds}
 \label{sec2}
  Let $M\subset \Pn {N}$ be a projective (compact) manifold of dimension $n$.
 This means that in the homogeneous coordinates of $\Pn N$, $M$ is given by the zeros
 of homogeneous polynomials. 
 In $M$ we consider the canonical orientation induced by its complex structure.
 Any Hermitian metric in $M$ induces such an orientation.
 In the following by $H^m(M,\Z)$ (resp. $H_m(M,\Z)$) we mean the image of the classical
 singular cohomology (resp. homology) in 
$H^m(M,\C)$ (resp. $H_m(M,\C)$). Therefore, we have killed the torsion.


 \subsection{De Rham cohomology}

 The de-Rham cohomology of $M$ is given by
 $$
 H^m_\dR(M):=
\frac{Z^m(M)}{dA^{m-1}(M)}
$$
where $A^{m}(M)$ (resp. $Z^{m}(M)$)  is the set of $C^\infty$ complex valued
differential $m$-forms (resp. closed  $m$-forms)
on $M$. From another side we have 
$$
H^m_\dR(M)=H^m(M,\C):=H^m(M,\Z)\otimes_\Z \C.
 $$
 We look $H^m(M,\Z)$ as a lattice in $H_\dR^m(M)$.
\subsection{Intersection form}
\label{intersectionform}
We have the Poincar\'e duality 
$$
P: H_m(M,\Z)\rightarrow H^{2n-m}(M,\Z), \ \int_{\delta'}P(\delta)=\langle \delta,\delta'\rangle
$$
under which for $n=m$ the non-degenerate bilinear intersection map $\langle \cdot ,\cdot \rangle$
in homology corresponds to the bilinear map
$$
\psi(\omega_1,\omega_2)=c\int_{M}\omega_1\wedge \omega_2, \
\omega_1,\omega_2\in H_\dR^m(M),
$$
where $c$ is a positive real number depending 
only on $M$ (note that any two Hermitian metric induce the same
orientation in $M$). It assures us that
for $\omega_1,\omega_2\in H^m(M,\Z)$ we have $c\int_{M}\omega_1\wedge \omega_2\in\Z$.
If we use the notion of 
singular $p$-chains (see \cite{gri}, p. 43) for the definition of integral on manifolds
we can assume that $c=1$.
\subsection{Hodge decomposition}
We have the Hodge decomposition
\begin{equation}
\label{11dez02}
H^m_\dR(M)=H^{m,0}\oplus H^{m-1,1}\oplus\cdots\oplus H^{1,m-1}\oplus H^{0,m},
\end{equation}
where $H^{p,q}\cong\frac{Z^{p,q}(M)}{dA^{p+q-1}(M)\cap Z^{p,q}(M)}$ and
 $Z^{p,q}_d$  is the set of $C^\infty$ 
 closed $(p,q)$-forms
on $M$. We have  the canonical inclusions
$H^{p,q}\rightarrow H^m_\dR(M)$ and one can  prove ~(\ref{11dez02}) using
harmonic forms, see M. Green's lectures \cite{gmv}, p. 14. We have the conjugation
mapping
$$
H_{\dR}^m(M)\rightarrow H_\dR^m(M),\ \omega\mapsto \bar\omega
$$
which leaves $H^m(M,\R)$ invariant and maps $H^{p,q}$ isomorphically to $H^{q,p}$.
By taking local charts it is easy to verify that

\begin{equation}
\label{151203}
\psi(H^{i,m-i},H^{m-j,j})=0 \hbox{ unless } i=j,
\end{equation}
\begin{equation}
\label{pirshodam}
(-1)^{i+\frac{m}{2}}\psi(a,\bar a)
>0,  \ \forall a\in H^{i,m-i},\ a\not=0.
\end{equation}
\subsection{Hodge conjecture}
One of the central conjectures in Hodge theory is the so called Hodge conjecture.
Let $m$ be an even natural number and $Z=\sum_{i=1}^s r_iZ_i$, where 
$Z_i,\ i=1,2,\ldots,s$ is a compact algebraic  subvariety of $M$ of
complex dimension $\frac{m}{2}$ and $r_i\in \Z$. 
 By Chow theorem every compact analytic subvariety of $M$ is  algebraic, i.e. it is given by zeros of polynomials (see \cite{gri}).
Using a
resolution map $\tilde Z_i\rightarrow M$, where $\tilde Z_i$ is a compact 
complex
manifold,  one can define an element
$\sum_{i=1}^s r_i[Z_i]\in H_m(M,\Z)$  which is called an algebraic cycle (see 
\cite{boha61}).
 Since the restriction to $Z$ of a $(p,q)$-form with
$p+q=m$ and $p\not=\frac{m}{2}$ is identically zero, an algebraic cycle
$\delta$ has the following property:
$$
\int_\delta H^{p,m-p}=0,\ \forall p\not = \frac{m}{2}.
$$
 A cycle $\delta\in H_m(M,\Z)$  with the above property
is called a Hodge cycle. The assertion of the Hodge conjecture is that for any 
Hodge cycle $\delta\in H_m(M,\Z)$ there is an integer $a\in \N$ such that 
$a\cdot \delta$ is an algebraic
cycle, i.e. there exist subvarieties $Z_i\subset M$ of dimension
$\frac{m}{2}$ and rational numbers $r_i$  such that $\delta=
\sum r_i[Z_i]$.  The difficulty of this conjecture
lies in constructing varieties just with their homological information.
The conjecture is false with $a=1$ (see \cite{gro69}).
In the literature one usually finds the notion of
a Hodge class which is an element in $H^{\frac{m}{2},\frac{m}{2}}\cap H^m(M,\Q)$. 
The Poincar\'e duality gives a bijection between the $\Q$-vector space 
generated by Hodge cycles and the $\Q$-vector space of Hodge
classes (see \cite{vo02} \S 11).
\section{The classifying space of polarized Hodge structure}
\label{sec3}
In this section we construct the classifying space of Hodge structures 
$D=D(m,h,\psi)$ with a fixed weight $m$, 
Hodge numbers $h=(h^{m-i,i},\ i=0,1,\ldots, m)$ and a polarization $\psi$ 
on a fixed freely generated $\Z$-module $\HH \Z {}$. 
One can look $D$ in two ways: First, as the space of Hodge filtrations 
 and second as a quotient $G/K$, where 
$G$ is a real Lie group and $K$ is a compact subgroup 
(not necessarily maximal). For a subring $\K$ of $\C$ we define
$$
\HH \K {}:=\HH \Z {}\otimes_\Z \K.
$$
\subsection{Polarized Hodge structures}
Let us be given a freely generated $\Z$-module $H(\Z)$. A Hodge structure of
weight $m$, $m\in{\mathbb N}$ and type\index{Hodge structure} 
$h=(h^{m,0},h^{m-1,1},\ldots,h^{0,m})\in {\mathbb N}^{m+1}$
on $H(\Z)$ is a decomposition
$$
H(\C):=H^{m,0}\oplus H^{m-1,1}\oplus\cdots\oplus H^{0,m}
$$
with ${\bar H}^{i,m-i}=H^{m-i,i}$ and  $h^{i,m-i}=dim_\C H^{i,m-i}$. We call
$h^{m-i,i}$ the Hodge numbers\index{Hodge numbers} 
and define $h^i:=h^{m,0}+h^{m-1,1}+\cdots, 
h^{m-i,i}$.  When $\HH{\Z}{}$ has a Hodge structure
of type $m$ we denote it by $\HH{\Z}{m}$. 

Instead of Hodge structures, we will use Hodge filtrations. The main reason
is that one can define de Rham cohomologies of algebraic varieties over a field $k$ 
and the associated 
Hodge filtrations in such a way that they coincide with the classical  notions in the case of 
$k=\C$ (see \cite{gro66}).
For $0\leq i\leq m$  we define
$$
F^i=F^i\HH\C m:=H^{m,0}\oplus H^{m-1,1}\oplus\cdots\oplus H^{i,m-i}.
$$
To recover the Hodge filtration we define
$H^{i,m-i}=F^i\cap \bar F^{m-i}$.


When we have a family of
Hodge structures parameterized by $\alpha \in I$ then we denote the 
Hodge structure associated to $\alpha\in I$ by $\HH{\alpha,\Z} m$.
This notation is also used replacing
$\Z$ with $\Q$, $\R$, $\C$ and so on. 
We write also  $F^i=F^i(\alpha,\C)=F^i(\alpha)$, $H^{m-i,i}=
\HH{\alpha,\C}{m-i,i}=\HH{\alpha}{i,m-i}$ and so on.  

A polarization \index{polarization}
$\psi=\langle \cdot ,\cdot \rangle $ for $H(\Z)$ is a non-degenerate bilinear map 
$H(\Z)\times H(\Z)
\rightarrow \Z$ symmetric if $m$ is even and skew if $m$ is odd
$$
\psi(a,b)=(-1)^m\psi(b,a), \ a,b\in \HH\Z m
$$ 
such that its complexification (we denote it again by $\psi$) 
satisfies (\ref{151203}) and (\ref{pirshodam}). All the data above is called 
a polarized Hodge structure of type $\Phi=(m,h,\psi)$.

\subsection{Hodge structure in cohomologies with real coefficients}
Every element 
$x\in \HH\R m$ can
be written in the form 
$$
x=x_{m,0}+x_{m-1,1}+\cdots+\overline{x_{m-1,1}}+
\overline{x_{m,0}}, \ x_{m-i,i}\in H^{m-i,i}
$$
with $\overline{x_{\frac{m}{2},
\frac{m}{2}}}=x_{\frac{m}{2},\frac{m}{2}}$ if $m$ is even. We define 
$H^i\subset \HH\R m, i<\frac{m}{2}$ to be the set of all 
$x_{m-i,i}+\overline{x_{m-i,i}},\ x\in \HH\C m$ and 
$H^{\frac{m}{2}}$, for $m$ even, to be the set of all $x_{\frac{m}{2},\frac{m}{2}},
 x\in \HH\C m$. 
These are  real subvector spaces of $\HH\R m$ and we have 
the following decomposition of $\HH\R m$:
\begin{equation}
\label{10122003}
\HH\R m=H^0\oplus H^1\oplus\cdots\oplus H^{\frac{m}{2}},
\end{equation}  
$$
x=(x_{m,0}+\overline{x_{m,0} })+(x_{m-1,1}+\overline{x_{m-1,1} })+
\cdots+x_{\frac{m}{2},\frac{m}{2}}.
$$
For $i<\frac{m}{2}$ 
the map $x_{m-i,i}\rightarrow x_{m-i,i}+\overline{  x_{m-i,i}}$ induces 
an isomorphism of $\R$-vector spaces $H^{m-i,i}$ and $H^i$. Multiplication
by $\sqrt{-1}$ in $H^{m-i,i}$ induces a map
\begin{equation}
\label{101203}
J_i:H^i\rightarrow H^i,\ J_i^2=-{\rm Id},
\end{equation}
$$
x=x_{m-i,i}+\bar x_{m-i,i},\ J_ix:=ix_{m_i,i}+\overline{ix_{m_i,i}}=i(x_{m-i,i}-\bar x_{m-i,i}).
$$
If $m$ is even then we define $J_{\frac{m}{2}}$ to be the identity. A $\C$-linear
map in $H^{m-i,i}$ corresponds to a $\R$-linear map in $H^i$ which
 commutes  with $J_i$.
\begin{prop}(Riemann relations)
\label{RiRe}
For $i,j\leq \frac{m}{2}$ we have: 
\begin{enumerate}
\item
 $\psi(H^i,H^j)=0$ unless $i=j$;
\item
$\psi(J_ix,J_iy)=\psi(x,y)$ for all $x,y\in H^i$; 
\item
If $m$ is odd then $(-1)^{\frac{m-1}{2}+i}\psi(x,J_ix)>0$ for all 
$x\in H^i,\ x\not=0$ ($\psi(x,x)=0$);
\item
If $m$ is even then $\psi(x,J_ix)=0$ and $(-1)^{\frac{m}{2}+i}\psi(x,x)>0$ for
all $x\in H^i,\ x\not=0$.
\end{enumerate}
\end{prop}
\begin{proof}
1. It is a direct consequence of (\ref{151203}). 2. Write $x=a+\overline a
, y=b+\overline b, a,b\in H^{m-i,i}$. 
$$
\psi(J_ix,J_iy)=\psi(\sqrt {-1}(a-\bar a),\sqrt{-1}(b-\bar b))=
-\psi(a-\bar a,b-\bar b)=
\psi(a,\bar b)+\psi(\bar a,b)=
$$
$$
\psi(a+\bar a, b+\bar b)=\psi(x,y).
$$
3,4. We use (\ref{151203}) to obtain 
$$
\psi(x,J_ix)=\psi(a+\bar a,\sqrt{-1}(a-\bar a))=\sqrt{-1}((-1)+(-1)^m)
\psi(a,\bar a),
$$
$$
\psi(x,x)=\psi(a+\bar a,a+\bar a)=\psi(a,\bar a)+\psi(\bar a,a)=(1+(-1)^m)
\psi(a,\bar a),
$$
and then use (\ref{pirshodam}).     
\end{proof}
Let $C:\HH\R m\rightarrow \HH\R m$ be defined in the following way
$$
Cx:=\left\{ \begin{array}{ll}
             (-1)^{\frac{m-1}{2}+i}J_ix &    m \hbox{ odd }, \\
             (-1)^{\frac{m}{2}+i}x  & m \hbox{ even},
             \end{array}
             \right.   x\in H^i.
$$
Note that $\psi(x,Cy)$ is a positive form on $\HH\R m$ (see \cite{de70}).
We will call the decomposition (\ref{10122003}) of $\HH\R m$ and the 
data (\ref{101203}) with the properties 1,2,3 and 4 of Proposition \ref{RiRe}
the polarized Hodge structure on $\HH\R m$. 
A polarized Hodge structure on $\HH\R m$  
gives  in a canonical way a polarized Hodge structure on $\HH\C m$.
To see this, for $i=\frac{m}{2}$ we define  $H^{\frac{m}{2},\frac{m}{2}}
=H^{\frac{m}{2}}\otimes
\C$ and for $i<\frac{m}{2}$ we define $H^{m-i,i}$ (resp. $H^{i,m-i}$) to be the
vector space generated by the eigenvectors of $J_i$ with the eigenvalue 
$\sqrt{-1} $(resp. $-\sqrt{-1})$. Now one can check 
(\ref{151203}) and (\ref{pirshodam}).


\begin{prop}
\label{14dez03}
For a polarized Hodge structure $\alpha$ on 
$\HH\Z{m}$, the set
\begin{equation}
\label{13.5.06}
V_{\alpha}:=
\{A\in GL(\HH\R m)\mid A \hbox{ respects the Hodge structure and } 
\end{equation}
$$
\psi(Ax,Ay)=\psi(x,y), \ AJ_i=J_iA, i=0,1,\ldots\}
$$
is a compact subgroup of $GL(\HH\R m)$.
\end{prop}
\begin{proof}
Define $P(x)=\psi(x,Cx), \ x\in \HH \R m$. 
We have $P(Ax)=P(x),\forall A\in 
V_{\alpha}$.
Therefore, elements of  $V_{\alpha}$ leaves the fibers $P^{-1}(c), c\in\R^+$ invariant. 
Since the fibers of $P$ are compact sets, this finishes the
proof of the proposition.
\end{proof}
\subsection{Generalized Jacobians}
Let us suppose that $m$ is odd and $i=\frac{m+1}{2}$.
 Then we have a canonical isomorphism $\HH\R{m}\rightarrow F^{\frac{m+1}{2}}$ of
$\R$-vector spaces.
Therefore, the projection $L$ of $\HH\Z{m}$ in $F^{\frac{m+1}{2}}$
is a lattice of
rank $2\dim_\C F^{\frac{m+1}{2}}$ and so we can define the compact complex  torus
$$
J_{\frac{m+1}{2}}\HH\Z{m}:=F^{\frac{m+1}{2}}/L.
$$
This is called the $\frac{m+1}{2}$-Jacobian variety \index{Jacobian variety}
 of $\HH\Z{m}$. 
In the case for which  $\HH\Z m$ is the integral cohomology of 
a smooth projective variety of dimension $m$,  
the torus $J_1\HH \Z 1, \ m=1$
(resp. $J_{n}\HH\Z {2n-1},\ m=2n-1$) is also  called the Albanese variety
\index{Albanese variety} 
(resp. Picard variety)\index{Picard variety}. 

\subsection{Griffiths domain}
We fix $m$, 
Hodge numbers $h=(h^{m-i,i},\ i=0,1,\ldots, m)$ and a non-degenerate bilinear map
$\psi=\langle \cdot,\cdot \rangle$ in $\HH \Z m$. The Griffiths domain $D$
 is the space of all 
 decompositions $\HH\C m:=H^{m,0}\oplus H^{m-1,1}\oplus\cdots\oplus H^{0,m}, \
 \dim H^{m-i,i}=h^{m-i,i}$ resulting a polarized Hodge structure on $\HH \Z m$. 
We define the compact dual $\check D$  of $D$ similar
to $D$ but without the condition (\ref{pirshodam}).
Note that $D$ is an open subset of $\check D$.
The compact dual $\check D$ of $D$ is an analytic variety
in the following way: First note that a Hodge structure on $\HH\Z m$ is
completely determined by the data:
\begin{equation}
\label{mehdi}
F^{[\frac{m}{2}]+1}=H^{m,0}\oplus H^{m-1,1}\oplus\cdots \oplus H^{[\frac{m}{2}]+1,
m-[\frac{m}{2}]-1},
\end{equation}
\begin{equation}
\label{ali}
F^{[\frac{m}{2}]+1}\cap\overline{ F^{[\frac{m}{2}]+1}}=\{0\},\  \psi(H^{m-i,i},
H^{m-j,j})=0,\ \forall i,j\leq m-[\frac{m}{2}]-1.
\end{equation}
($H^{i,m-i}$ is defined to be $\overline{H^{m-i,i}}$ and 
in the case $m$ even $H^{\frac{m}{2},\frac{m}{2}}$ is the $\psi$-orthogonal complement of $F^{[\frac{m}{2}]+1}+\overline{ F^{[\frac{m}{2}]+1}}$). 
The decomposition (\ref{mehdi}) determines a point in the 
complex grassmannian manifold 
$\gr:=\gr (h^{m,0},\HH \C m)\times\cdots\times \gr (h^{[\frac{m}{2}]+1,
m-[\frac{m}{2}]-1},\HH \C m)$. The first condition in (\ref{ali}) determines
an open subset of $\gr$ and the second condition an analytic 
subvariety of
$\gr$. 
\subsection{$D$ as a quotient of  real Lie groups}
For a subring $\K$ of $\R$ define
$$
\Gamma_\K:=\Aut(H(\K),\langle \cdot,\cdot\rangle):=
$$
$$
\{A:H(\K)\to H(\K) \mid A \hbox{ is $\K$-linear and } \forall x,y\in H(\K), \ \langle Ax,Ay\rangle=\langle x,y\rangle \}.
$$
The group $\Gamma_\R$ acts from the left on $D$ in a canonical way. 
For a fixed  point $p_0\in D$ define
$$
V=V_{p_0}:=\{a\in \Gamma_\R\mid a.p_0=p_0 \}.
$$ 
According to Proposition \ref{14dez03}, $V$ is a compact
subgroup of $\Gamma_\R$.
The map 
$$
\alpha: \Gamma_\R/V\rightarrow D,\ \alpha(a)=a\cdot p_0
$$ is 
an isomorphism and so we may identify $D$ with $\Gamma_\R/V$.
Note that $V$ may not be maximal. 

Since $\Gamma_\Z$ is discrete in $\Gamma_\R$, i.e. 
it has the discrete
topology as a subset of $\Gamma_\R$, 
the group $\Gamma_\Z$ acts discontinuously on $D$, i.e.
for any two compact subset $K_1,K_2$ in $D$ the set 
$$
\{A\in\Gamma_\Z\mid A(K_1)\cap K_2\neq \emptyset\} 
$$ 
is finite.
%
The set
 $\Gamma_\Z\backslash D$ is the moduli of polarized Hodge structures of
 type $\Phi=(m,h,\langle\cdot,\cdot\rangle)$ and it has a
canonical structure of a complex analytic space. 
\subsection{Tangent space of $D$}
Let
$$
\lieg_\K=
{\mathrm Lie}(\Gamma_\K):=
\{N\in \End_\K(\HH \K m)  \mid \langle Nx,y\rangle +\langle x,Ny\rangle =0,\ 
\forall x,y\in \HH \K m\}.
$$
To each $\alpha\in D$, 
there is a natural filtration in $\lieg_\C$ 
$$
F^{i}\lieg_\C=\{N\in\lieg_\C\mid N(F^p)\subset F^{p+i},\ \forall p\in\Z\},
\ i=0,-1,-2,\ldots,
$$
where $F^\bullet$ is the Hodge filtration associated to $\alpha$. 
We get a natural filtration of the tangent bundle
$$
T^h_\alpha D:=F^{-1}(\lieg_\C)/ F^0(\lieg_\C) 
\subset F^{-2}(\lieg_\C)/ F^0(\lieg_\C)
\subset\cdots\subset \lieg_\C/ F^0(\lieg_\C)=T_\alpha D.
$$
One usually calls  $T^h_\alpha D$ the horizontal tangent bundle. 
\index{horizontal tangent bundle}
When $m=1$ then the horizontal and usual
tangent bundles are the same.

The subgroup $V$ is connected and  is contained in a unique maximal compact subgroup 
$K$ of $\Gamma_\R$. 
When $K\neq V$, then there is a fibration of $D=\Gamma_\R/V\rightarrow \Gamma_\R/K$
with compact fibers isomorphic to $K/V$ which are complex subvarieties of $D$. In this case
we have  $T_\alpha(D)=T^h_\alpha (D)\oplus T^v_\alpha(D)$, where $T^v(D)$ restricted to a
 fibre of $\pi$ coincides with the tangent bundle of that fibre.
For more information see \cite{caka78} and the references there.

\subsection{Few good cases}
\label{fewcases}
Here is the classification of all cases in which $V$ is a maximal compact subgroup of $\Gamma_\R$,
i.e. it is compact and there is no other compact subgroup of $\Gamma_\R$ containing $V$.  
\begin{prop}
(\cite{kaus04} page 4 or \cite{caka77})
The group $V$ is maximal in $\Gamma_\R$ only in  one of the following cases:
\begin{enumerate}
\item
$m=2a+1, h^{p,n-p}=0$ if $p\not= a,a+1$;
\item
$m=2a, h^{a+1,a-1}\leq 1$ and $h^{p,m-p}=0$ if $p\not=a-1,a,a+1$.
\end{enumerate} 
\end{prop} 
In the first case $D$ is the Siegel domain (see \cite{fr83, kl90}). 

A subgroup of $\Gamma_\R$  is called arithmetic if it
is commensurable with $\Gamma_\Z$. Here two subgroups A and B of a group are commensurable when their intersection has finite index in each of them. 
For the case in which $V$ is maximal in $\Gamma_\R$, the quotient
$D:=\Gamma_\R/V$ is a Hermitian symmetric domain
(see for instance \cite{mil03}). 
For an arithmetic subgroup $\Gamma$ of $\Gamma_\R$, the compactification of
$\Gamma\backslash D$ is done by I. Satake, A. Borel and W. Baily 
(\cite{babo66, sa60}).
In fact the compactification due to Borel and Baily gives us an algebraic 
structure on $\Gamma\backslash D$.
For the case in which $V$ is not maximal, partial compactifications are done 
by E. Cattani, A. Kaplan, A. Ash, D. Mumford, M. Rapoport, Y. S. Tai, K. Kato and S. Usui (see the references of \cite{kaus00, kaus04}).
\section{Period map}
\label{sec4}
Roughly speaking the period map associate to each variety its polarized Hodge structure
and hence a point in the Griffiths domain.
\subsection{Ehresmann fibration theorem}
For a family of  algebraic varieties depending on a parameter one can always find a Zariski 
open subset $U$ in the parameter space in such a way that the varieties with corresponding
parameter in $U$ are topologically the same (see for instance \cite{ve76}, corollary 5.1).  
However, they may have different analytic structures and Hodge structures.
\begin{theo}
(Ehresmann's Fibration Theorem \cite{eh47}). 
Let $f:Y\rightarrow B$ be a proper submersion between the manifolds $Y$ and $B$.
 Then $f$ fibers $Y$ locally trivially i.e., for every point $b\in B$ 
 there is a neighborhood $U$ of $b$ and a $C^\infty$-diffeomorphism 
 $\phi :U\times f^{-1}(b) \rightarrow f^{-1}(U)$ such that 
 $f\circ\phi =\pi _1 =$ the first projection. Moreover if $N\subset Y$ 
 is a closed submanifold such that $f\mid _N$ is still a submersion then 
 $f$ fibers $Y$ locally trivially over N i.e., the diffeomorphism  
 $\phi$ above  can be chosen to carry  $U\times (f^{-1}(b)\cap N)$ 
 onto $f^{-1}(U)\cap N$.
\end{theo}
The map $\phi$ is called the fiber bundle trivialization map. 
Ehresmann's theorem can be rewritten for manifolds with boundary and also 
for stratified analytic sets. In the last case the result is due to 
R. Thom, J. Mather and J. L.  Verdier (see \cite{mat71, ve76}).

Let  $\lambda$ be a path in $B$ connecting $b_0$ to $b_1$ and
defined
up to homotopy. 
Ehresmann's Theorem gives us a unique map 
$h_\lambda: f^{-1}(b_0)\rightarrow f^{-1}(b_1)$  defined up to homotopy. In particular, 
for $b:=b_0=b_1$ 
we have  the action of $\pi_1(Y,b)$ on the homology group $H_n(f^{-1}(b),\Z)$. 
The image of $\pi_1(Y,b)$ in $\Aut (H_n(f^{-1}(b),\Z))$ is called the monodromy group.

\subsection{period map}
\label{darhavapeima2006}
Let us be given a holomorphic map between projective (compact) 
varieties  $f:X\to S$. Let $S'$ be the locus of points $t\in S$ such that $f$
is not a submersion. According to Ehresmann's theorem the mapping
$f$ restricted to $T:=S\backslash S'$ is a $C^\infty$ fiber bundle. 
Fix a point $t_0\in T$ and identify $(H^m(f^{-1}(t_0),\Z), \langle \cdot ,\cdot \rangle)$ 
with the polarized $\Z$-module in the definition of $D$. Let
$$
\Gamma=\Im(\pi_1(T,t_0)\mapsto \Gamma_\Z)
$$
be the monodromy group associated to $f$.
We have the well-defined period map 
$$
\per: T\rightarrow \Gamma \backslash D.
$$
I believe that the monodromy group $\Gamma$ for a complete $f$ (see \S\ref{kst}) is arithmetic but I do not know any reference or
proof for this fact. To avoid this problem we may use $\Gamma_\Z\backslash D$ instead of 
$\Gamma\backslash D$. 
\subsection{Griffiths transversality theorem}
The Griffiths transversality theorem is originally stated using the Gauss-Manin
connection on the $m$-th cohomology bundle on $T$. It implies that  
the image of the derivation of the period map is in the horizontal tangent space 
$T^h_{\per(t)}D$ 
of $D$.
\section{Automorphic functions}
\label{sec5}
Automorphic functions, from the point of view of this text, are connecting
functions between coefficient spaces of algebraic varieties and the corresponding 
space of integrals/Hodge structures. 
\label{autfun}
\subsection{Positive line bundles and Kodaira embedding theorem}
A line bundle $L$ on a compact complex manifold $A$ is called negative if the zero
section $A_0$ of $L$ can be contracted to a point in the context of analytic varieties, i.e
there is an analytic map from a neighborhood of $A_0$ in $L$ to a singularity $(X,0)$ such that
it is a biholomorphism outside $A_0$ and the inverse image of $0$ is $A_0$.
Naturally, a line bundle is called positive if its dual is negative.  
\begin{theo}
Let $L$ be a positive line bundle on a complex compact manifold $A$. Then there is an
integer $n$ and global holomorphic sections $s_0,s_1,\ldots,s_N$ of $L^n$ such that
the mapping
$$
A\to \Pn N,\ x\mapsto [s_0(x);s_1(x);\cdots:s_N(x)]
$$
is an embedding.
\end{theo}
As far as I know this is the only way for giving an algebraic structure
to a complex manifold. For further information the reader is referred to \cite{gra62, camo}.
\subsection{Automorphy factors}
Let $D$ be a complex manifold and $\Gamma$ be a subgroup of biholomorphism group of 
of $D$. A holomorphic automorphy factor for $\Gamma$ on $D$ is a map
$\ja :D\times \Gamma\rightarrow \C^*$  which for fixed $A\in \Gamma$ is holomorphic in $x\in D$
and which satisfies the identity
$$
\ja (x, A\cdot B)=\ja(x,A)\ja(A\cdot x, B),\ \forall A,B\in \Gamma, \ x\in D.
$$
Two automorphic factors  $\ja_1$ and $\ja_2$ are called equivalent  if
there is a group homomorphism $a:\Gamma\to (\C^*,\cdot)$ such that $\ja_1(x,\cdot)=a(\cdot)
\ja_2(x,\cdot)$
for all $x\in D$. 
If both $D$ and $D':=\Gamma\backslash D$ are smooth varieties then the equivalence class  
of an automorphy factor $\ja$ corresponds to a  line bundle $L_\ja$ on $D'$ . Conversely, every line
bundle on $D'$ is obtained by an automorphy factor.  A global holomorphic section  of 
$L_{\ja}^n$ corresponds to a holomorphic function $s$ in $D$ such that
\begin{equation}
\label{4.5}
s(A\cdot x)=\ja(x,A)^ns(x), \forall \  x\in D,\ A\in \Gamma.
\end{equation}
In an arbitrary case in which $D'$ may not be a complex manifold it is natural to say that
a holomorphic function on $D$ is an automorphic function of weight $n$ 
if it satisfies (\ref{4.5}) and a certain growth condition (depending on the situation).
Usually, for a holomorphic function $s$ in $D$ one defines the slash operator
$$
s|_nA(\cdot ):=\ja(\cdot,A)^{-n}s(\cdot),\ A\in\Gamma
$$
which satisfies
$$
s|_n(AB)=(s|_nA)|B,\ A,B\in\Gamma
$$
and one rewrites (\ref{4.5}) in the form
\begin{equation}
\label{4.5.6}
s|_nA=s,\ \forall A\in \Gamma.
\end{equation}
If we wish  to find a canonical embedding of $D'$ in some projective space and we want to use 
the idea behind Kodaira embedding theorem then we have to find an automorphy factor
$\ja$, a positive integer $n$ and automorphic functions $s_0,s_1,\ldots,s_N$ of weight $n$ such that
the map
$$
D\to \Pn N,\ x\mapsto [s_0(x);s_1(x);\cdots:s_N(x)] 
$$  
is an  embedding.

If $D$ is a domain in $\C^n$ then the determinant of the Jacobian of $h\in \Gamma$ at $x$, denote
it by $\det\jacob(x,h)$ is an automorphy factor. The corresponding line bundle in 
$D'$ is the canonical line bundle in $D'$, i.e. the wedge product of
the cotangent bundle of $D'$ $n$-times. 
\subsection{Poincar\'e  series}
Let us consider a holomorphic function $s$ in $D$ which may not satisfy the equality
(\ref{4.5}). To $s$ one can associate the formal series
$$
\tilde s(x):=\sum_{A\in\Gamma} s|_nA.
$$ 
It satisfies the property (\ref{4.5.6}) but it may not be convergent.
The Borel-Baily theorem says that for $D$ a Hermitian symmetric domains
and for suitable $s$ the above series is convergent.
\subsection{Borel-Baily Theorem}
Let ${\rm Hol}(D)^+$ be the connected component of the group if biholomorphisms of a variety $D$ containg the identity. 
\begin{theo}
(Borel-Baily)
Let $\Gamma\backslash D$ be the quotient of a Hermitian symmetric
domain by a torsion free arithmetic subgroup of ${\rm Hol}(D)^+$ then
there is a positive integer $n$ such that the automorphic forms of weight
$n$ gives us an embedding of $\Gamma\backslash D$ in some projective space.
\end{theo}
The automorphy factor in the above theorem comes, as usual, from the canonical
bundle of $\Gamma\backslash D$. The proof is mainly based on the study of
automorphic functions and convergence of the Poincar\'e series (see \cite{babo66}, \S 5).

Shimura Varieties are special cases of quotients $\Gamma\backslash D$. They represent
certain moduli spaces in Algebraic geometry. For further information on this subject
the reader is referred to \cite{mil03}. 
\subsection{Final note on the moduli of polarized Hodge structures}
Let $D$ be the Griffiths domain and $\Gamma$ be an arithmetic subgroup 
of $\Gamma_\R$. We have the following vector bundle on $D':=\Gamma\backslash D$:
$$
H^m:= \cup_{\alpha\in D'}H^m(\alpha,\C).
$$ 
For $0\leq i\leq m$ it has the subbundle  $F^i:=\cup_{\alpha\in D'}F^i(\alpha,\C)$.
The wedge product of $F^i$, $\rank(F^i)$ times, gives us a line bundle in $D'$ and
hence an automorphic factor in $D$. Except in the few cases mentioned in \S \ref{fewcases}
these line bundles are not positive, in the sense that they have not enough holomorphic
sections in order to embed $D'$ in some projective space. 
\section{A new point of view}
\label{sec6}
In the construction of the Griffiths domain $D$, we have considered many Hodge structures
that may not come from geometry. In the simplest way we may define that a Hodge structure
comes from geometry  if it arises in the $m$-th cohomology of some smooth algebraic variety.
However, A. Grothendieck in \cite{gr70} p. 260 gives an example in which a Hodge
structure comes in a certain way from geometry but it is not included in our
premature definition. In the case of Hodge structures arising from Riemann surfaces
of genus $g\geq 2$ the Griffiths domain is of dimension $\frac{g(g+1)}{2}$  and its subspace consisting of Hodge structures coming from Riemann surfaces is of dimension $3g-3$.
The conclusion is that it would not be a reasonable idea to look for certain 
algebraic structures for $\Gamma\backslash D$.  Instead, we propose a 
point of view which is explained in this section. For simplicity, we explain it in the case of 
hypersurfaces. 
\subsection{Kodaira-Spencer Theorem on deformation of hypersurfaces}
\label{kst}
For a given smooth hypersurface $M$  of degree $d$ in $\Pn {n+1}$
is there  any deformation of $M$ which is not embedded in $\Pn {n+1}$? 
 The
answer to our question is no, except for some few cases.  
It is given by Kodaira-Spencer Theorem which we are going to explain
it in this section. For the proof and more information on deformation of complex manifolds the
reader is referred to \cite{kod86}, Chapter 5.

Let $M$ be a complex manifold and $M_t, t\in B:=(\C^s,0), \ M_0=M$ be
a deformation of $M_0$ which is topologically trivial over $B$. 
We say that the parameter space $B$ is effective if
the Kodaira-Spencer map
$$
\rho_0:T_0B\rightarrow H^1(M,\Theta)
$$ 
is injective, where $\Theta$ is the sheaf of vector fields on $M$. It is
called complete if other families are obtained from $M_t, t\in B$
in a canonical way (see \cite{kod86}, p. 228).
\begin{theo}
If $\rho_0$ is surjective at $0$ then $M_t, t\in B$ is complete.
\end{theo}
If one finds an effective deformation
of $M$ with $\dim B=\dim_\C H^1(M,\Theta)$ then $\rho_0$ is surjective and so
by the above theorem it is complete.  

Let us now $M$ be a smooth hypersurface of degree $d$ in the projective
space $\Pn {n+1}$. 
Let $T$ be the projectivization
of the  coefficient space of smooth hypersurfaces in $\Pn {n+1}$.
In  the definition of $M$ one has already
$\dim T =\binom{n+1+d}{d}-1$ parameters, from which  only
$$
m:=\binom{n+1+d}{d}-(n+2)^2 
$$ 
are not obtained by linear transformations of $\Pn {n+1}$. 
\begin{theo}
Assume that $n\geq 2$, $d\geq 3$ and $(n,d)\not =(2,4)$. Then   
there exists a germ at $M$ of $m$-dimensional ($m=\dim_\C H^1(M,\Theta)$) smooth
subvariety  of $T$, say $B$,   such
that the Kodaira-Spencer map on $B$  is injective and so the corresponding deformation 
is complete.      
\end{theo}
For the proof see \cite{kod86} p. 234.  
Let us now discuss the exceptional cases. For $(n,d)=(2,4)$ we have
$19$ effective parameter but $\dim H^1(M,\Theta)=20$. The difference comes from
a non algebraic deformation of $M$ (see \cite{kod86} p. 247).  In this case $M$ is
a $K3$ surface.
For $n=1$,  we are talking about the deformation theory of a Riemann surface.
According to Riemann's well-known formula, the complex structure of 
a Riemann surface of genus $g\geq 2 $ depends on $3g-3$ parameters which
is again $\dim  H^1(M,\Theta)$ (\cite{kod86} p. 226).
\subsection{Tame polynomials}
Let $\alpha=(\alpha_1,\alpha_2,\ldots,\alpha_{n+1})\in\N^{n+1}$  
and assume that
the greatest common divisor of all $\alpha_i$'s is one.
We consider  a parameter ring $\ring:=\C(t),\ t=(t_1,t_2,\ldots,t_s)$. We also consider
the polynomial ring  
$\ring[x]:=\ring[x_1,x_2,\ldots,x_{n+1}]$ as a graded
algebra with $deg(x_i)=\alpha_i$.  A polynomial $f\in\ring[x]$ is called
 a quasi-homogeneous polynomial of degree $d$ with respect to the grading
$\alpha$ if
$f$ is a linear combination of monomials of the type
$x^\beta:=x_1^{\beta_1}x_2^{\beta_2}\cdots
x_{n+1}^{\beta_{n+1}},\alpha.\beta:=\sum_{i=1}^{n+1}\alpha_i\beta_i=d$. For
an arbitrary polynomial $f\in\ring[x]$ one can write in a unique way
$f=\sum_{i=0}^df_i,\ f_d\not=0$,
where $f_i$ is a quasi-homogeneous polynomial of degree
$i$. The number $d$  is called the degree of
$f$. 

Let us be given a polynomial $f\in\ring[x]$.
We assume that $f$ is a tame polynomial. In this text this means
that  there exist natural numbers $\alpha_1,\alpha_2,\ldots,\alpha_{n+1}\in\N$ such that
the Milnor vector space
$$
V_g:=\ring[x]/<\frac{\partial g}{\partial
  x_i}\mid i=1,2,\ldots,n+1> 
$$
is a finite dimensional $\ring$-vector space,
where $g=f_d$ is the last quasi-homogeneous piece
of $f$ in the graded algebra $\ring[x],\ {\mathrm deg}(x_i)=\alpha_i$.

We choose a basis
$x^I=\{x^\beta\mid\beta\in I\}$ of monomials for the Milnor $\ring$-vector space
and define the Gelfand-Leray $n$-forms
\begin{equation}
\label{21may2004}
\omega_\beta:=\frac{x^\beta dx}{df}, \ \beta\in I,
\end{equation}
where $dx=dx_1\wedge\cdots \wedge dx_n\wedge dx_{n+1}$ (see \cite{arn}). 
Let $\as_0=\C^s$, $\Delta\in \C(t)$ be the discriminant of $f$ and $T:=\as_0\backslash 
({\rm Zero}(\Delta)\cup {\rm Pole}(\Delta))$.    It turns out
that  $L_t:=\{f_t=0\},\ t\in T$ is a topologically trivial family of affine 
hypersurfaces, 
 where $f_t$ is obtained from $f$ by fixing
the value of $t$.  The differential forms 
$\omega_\beta,\beta\in I$ defined in (\ref{21may2004}) restricted to $L_t,\ t\in T$ form a basis
of $H_\dR^n(L_t)$.
The reader is referred to \cite{ho06-1, mo} for all unproved statements in this section.
In these articles we have also given algorithms which calculate a basis of the de Rham 
cohomology compatible with the mixed Hodge structure of the affine variety 
$L_t$. Mixed Hodge structures generalize the classical Hodge structures for arbitrary
varieties which may be non compact and singular. The reader is referred to \cite{de71, de74} on
this subject.   

The reader may have already noticed that the main reason for us to use arbitrary weights $\alpha_i$, 
is to put the example (\ref{family}) and the hypersurfaces discussed in 
 \S \ref{kst} into one context.
\subsection{Poincar\'e series}
Despite the fact that the differential forms $\omega_\beta,\ \beta\in I$ may not have any compatibility
with the mixed Hodge structure of $L_t$, we may still ask for the convergence of
Poincar\'e type series explained bellow. Define the period matrix
$$
\per(t)=
\begin{pmatrix}
\int_{\delta_1} \omega_1 & \int_{\delta_1} \omega_2  
 & \cdots &  \int_{\delta_1} \omega_\mu  \\
\int_{\delta_2} \omega_1 & \int_{\delta_2} \omega_2 
   & \cdots &  \int_{\delta_2} \omega_\mu  \\
\vdots  & \vdots  & \vdots  & \vdots  \\ 
\int_{\delta_\mu} \omega_1 & \int_{\delta_\mu} \omega_2  
  & \cdots &  \int_{\delta_\mu} \omega_\mu
\end{pmatrix},
$$
where $\delta=(\delta_{i},\ i=1,2,\ldots,\mu)$ is a basis of the freely generated $\Z$-module 
$H_n(L_t,\Z)$ and $\{\omega_i\mid i=1,2,\ldots,\mu\}=\{\omega_\beta\mid \beta\in I\}$.
Assume that the intersection matrix in the basis $\delta$ is $\Psi$. Different choices of
the basis $\delta$ will lead to the action of
$$
\Gamma_\Z:=\{A\in\Mat(\mu\times \mu, \Z)\mid A\Psi A^\tr=\Psi\}
$$
from the left on $\per(t)$. For a meromorphic function
$P:\Mat(\mu\times\mu, \C)\rightarrow \C$ define
$$
\Gamma_P=\{A\in \Gamma_\Z \mid P(Ax)=P(x),\ \forall x\in \Mat(\mu\times\mu, \C) \}.
$$
The Poincar\'e series of $P$ in the context of this section is
defined to be
$$
\check {P}(t):=\sum_{A\in\Gamma_P\backslash\Gamma_\Z }P(A\cdot \per(t)).
$$
If $\check P$ is convergent then by definition it is
one valued in $T$. The main question we pose
in this text is the convergence of $\check P$ and the description of the
sub algebra of $\C(t)$ generated by those convergent $\check P$ which extend meromorphically to $\as_0$. 
For the discussion of these problems in the case
\begin{equation}
\label{khodaya}
f=y^2-4t_0(x-t_1)^3+t_2(x-t_1)+t_3,\ t\in\C^4
\end{equation}
the reader is referred to \cite{ho06-2,ho06-3}.
\subsection{A convergence criterion}
In this section we describe a method which is used to prove that
certain Poincar\'e series are convergent. Similar methods can be found 
 in 
\cite{fr83} Satz 4.3, \cite{fr90} Lemma 5.1 p. 55 and 
\cite{babo66} Theorem 5.3 p. 49. 
In the third reference the authors associate their
convergence theorem  to Harish-Chandra.

Let $M$ be a complex manifold and $ds^2$ be a Hermitian-K\"ahlerian metric in $M$ and
let $dz$ be the associated volume form.  
Let $K\subset U\subset M$, where $K$ is a compact set and $U$ is
an open set. There is a real positive constant $C$ depending only on $(M, ds^2)$ and $K$ 
 such that for all holomorphic functions 
$f:M\to \C$  we have:
$$
|f(a)|^2<C\int_U |f(z)|^2dz, \ \forall a\in K,
$$
(see \cite{fr83} Satz 4.3, Hilfsatz 2, \cite{fr90} Lemma 5.1). Note that every Hermitian-K\"ahlerian form
can be written in a local holomorphic 
coordinates system  $(z_1,z_2,\ldots,z_n)$ as $dz_1\otimes d\bar z_1+\cdots+
dz_n\otimes d\bar z_n$.

Let $\Gamma\subset \Aut(M)$ such that $\Gamma$ acts discontinuously 
on $M$. This implies that for an arbitrary $a\in M$ the stabilizer  
$\Gamma_a:=\{A\in \Gamma\mid Aa=a\}$ of $a$ is finite and 
there exists an open neighborhood $U$ of $a$ such that 
$$
A\in\Gamma,\ A(U)\cap U\not =\emptyset\Rightarrow A\in\Gamma_a,\ 
$$
$$
A\in \Gamma_a\Rightarrow A(U)=U.
$$
We assume that the Hermitian metric of $M$ 
is invariant under the action of $\Gamma$. 
Let us take a holomorphic function $f$ on $M$.
We claim that if 
$$
\int_{M}|f|^2dz<\infty
$$
then the  Poincar\'e series 
$
\check f(z)=\sum_{A\in\Gamma}f(Az)
$
converges uniformly in $z$. This follows from the equalities:
$$
\sum_{A\in\Gamma}|f(Az)|^2\leq C\sum_{A\in \Gamma}\int_U|f(Az)|^2dz=
C\sum_{A\in\Gamma} \int_{A^{-1}U}|f(z)|^2dz\leq C\cdot  \#\Gamma_a \int_{M}|f|^2dz<\infty.
$$

\subsection{References for further investigation}
In this section I give a  list of articles and books which may be useful for
further development of the ideas explained in the present text. Of course
the reader will find much more literature, if he/she looks for the reference
citation or review citations of the mentioned works in Mathematical Review 
or Zentralblatt Mathematik. 

For the literature on arithmetic and algebraic groups the book \cite{bo69} 
is a good source of information. We have also the book  \cite{sp98} on algebraic
groups. The original paper of Baily and Borel \cite{babo66} can be served as a source
for Poincar\'e/Eisenstein  series on Hermitian symmetric domains. 
Some simplifications
are done in \cite{cas97}.   
An explicit construction of 
resolutions of the Borel-Baily compactification is given in \cite{amrt, kkms}. 
Further developments in the compactification problem is sketched in \cite{boji02}.
For the arithmetic point of view for the quotients  of Hermitian symmetric domains
by arithmetic groups the reader is referred to the original papers of Shimura \cite{shi70, shi1970} and
Deligne's papers  \cite{de70, de77}. For the study of the cohomology of
Shimura varieties one can mention the article \cite{hazu94} and two other papers with the
same title. The text \cite{mil03} can be served as an up-to-date  exposition
of the subject.

  In the Hodge theory side of the subject the expository article of Griffiths \cite{gr70}
is still a good source of information. See also \cite{gr69} for Deligne's report on
Griffiths works.  The compactification problem in Hodge theory can be seen as
the determining the limit of Hodge structures. For this problem see 
\cite{stzu85, sch73, caka77,caka78, kas85, cakasc86}. Recently, there have been attempts to look
at
the compcatification problem from log-geometry point of view (see \cite{kaus00, kaus04}). 
For the literature on log geometry the reader is referred to \cite{ka88, kana99,  kmn02}


\def\cprime{$'$} \def\cprime{$'$} \def\cprime{$'$}

\bibliographystyle{plain}
\end{document}